\documentclass[12pt]{amsart}
\usepackage{graphicx,mathrsfs,color,enumerate}
\usepackage[colorlinks=true]{hyperref}

\definecolor{mycolor}{rgb}{0.,0.0,1}

\hypersetup{colorlinks=true,linkcolor=mycolor,citecolor=mycolor,%
            urlcolor=mycolor,pdfpagemode=false,letterpaper=true,%
            pdfauthor={Mohamad A. Hindawi},%
            pdftitle={On the Filling Invariants at Infinity of Hadamard Manifolds}}

\allowdisplaybreaks[1]%
\hypersetup{pdfstartview=FitH}

\newtheorem{thm}{Theorem}[section]
\newtheorem{cor}[thm]{Corollary}
\newtheorem{lem}[thm]{Lemma}
\newtheorem{prop}[thm]{Proposition}
\theoremstyle{definition}
\newtheorem{defn}[thm]{Definition}
\theoremstyle{remark}
\newtheorem{rem}[thm]{Remark}
\numberwithin{equation}{section}
\newcommand{\norm}[1]{\left\Vert#1\right\Vert}
\newcommand{\abs}[1]{\left\vert#1\right\vert}

\newcommand{\Real}{\mathbb R}

\newcommand{\To}{\longrightarrow}

\newcommand{\gl}{\mathfrak{g}}
\newcommand{\kl}{\mathfrak{k}}
\newcommand{\pl}{\mathfrak{p}}
\newcommand{\al}{\mathfrak{a}}

\DeclareMathOperator{\vol}{vol}
\newcommand{\divv}{\operatorname{div}}
\newcommand{\inj}{\operatorname{inj}}
\newcommand{\sn}{\operatorname{sn}}
\newcommand{\cs}{\operatorname{cs}}
\newcommand{\ct}{\operatorname{ct}}

\begin{document}

\title[On the Filling Invariants at Infinity of Hadamard Manifolds]
{On the Filling Invariants at Infinity\\ of Hadamard Manifolds}%
\author{Mohamad~A. Hindawi}%
\address{Department of Mathematics\\University of Pennsylvania\\Philadelphia, PA 19104-6395\\USA}%
\email{\href{mailto:mhindawi@math.upenn.edu}{mhindawi@math.upenn.edu}}%

\thanks{}%

\begin{abstract}
We study the filling invariants at infinity $\divv_k$ for Hadamard
manifolds defined by Brady and Farb in~\cite{Brady98FIAIFMONC}.
Among other results, we give a positive answer to the question
they posed: whether these invariants can be used to detect the
rank of a symmetric space of noncompact type.
\end{abstract}
\maketitle

\section{Introduction} \label{S:Introduction}
Asymptotic invariants have been used perviously as basic tools to
study large scale geometry of various spaces. Precedents include
Gromov in~\cite{Gromov93AIOIG}, and Block and Weinberger
in~\cite{Block97LSHTAG}. In~\cite{Knieper97OTAGONCM} Knieper,
generalizing the work of Margulis in the case of negative
curvature, studied the growth function of the volume of distance
spheres in Hadamard manifolds and found an unexpected relation to
the rank of the manifold.

In~\cite{Gersten94QDOGICS} Gersten studied the divergence of
geodesics in ${\rm CAT}(0)$ spaces, and gave an example of a
finite ${\rm CAT}(0)$ \text{$2$-complex} whose universal cover
possesses two geodesic rays which diverge quadratically and such
that no pair of geodesics diverges faster than quadratically.
Adopting a definition of divergence which is a quasi-isometry
invariant, Gersten introduced a new invariant for geodesic metric
spaces, which we refer to as $\divv_0$. In~\cite{Gersten94DI3MG} Gersten
used the $\divv_0$ invariant to distinguish the quasi isometry type of
graph manifolds among all closed Haken \text{$3$-m}anifolds.

Using the same trick introduced by Gersten
in~\cite{Gersten94QDOGICS} to get a quasi-isometry invariant,
Brady and Farb in~\cite{Brady98FIAIFMONC} introduced a family of
new quasi-isometry invariants $\divv_k(X)$ for $0 \leq k \le n-2$,
for a Hadamard manifold $X^n$. These invariants were meant to be a
finer measure of the spread of geodesics in $X$.

The precise definition will be given below, but roughly the
definition of $\divv_k(X)$ is as follows: Find the minimum volume
of a \text{$(k+1)$-ball} which is needed to fill a
\text{$k$-sphere} which lies on the \text{$(n-1)$-dimensional}
distance sphere $S_{r}(x_0)$ in $X$. The filling is required to
lie outside the \text{$n$-dimensional} open ball
$B_{r}^{\circ}(x_0)$. The $\divv_k(X)$ invariant measures the
asymptotic behavior of the volume of the filling as $r \to \infty$
when the volume of the \text{$k$-sphere} grows polynomially in
$r$.

After computing some of these invariants for certain Hadamard
manifolds, Brady and Farb posed the following two questions.
\hyperdef{question}{1}{}{\\\bf Question~1.} Can the $\divv_k(X)$
invariants be used to detect the rank of a noncompact symmetric
space $X$?\\
{\bf Question~2.} What symmetric spaces can be distinguished by
the invariants $\divv_k$?

In this paper we study the $\divv_k$ invariants for various
Hadamard manifolds including symmetric spaces of noncompact type.
The first result we obtain is the following theorem.

\begin{thm} \label{T:higherRankSymmetricSpaces}
    If $X$ is an $n$-dimensional symmetric space of nonpositive
    curvature and of rank $l$, then $\divv_k(X)$ has a  polynomial growth of
    degree at most $k+1$ for every $k \geq l$.
\end{thm}

Brady and Farb showed in~\cite{Brady98FIAIFMONC} that
$\divv_{k-1}(X)$ is exponential when $X = H^{m_1} \times \dots
\times H^{m_k}$ is a product of $k$ hyperbolic spaces. The idea
was to show that there are quasi-isometric embeddings of $H^{(m_1
+ \dots + m_k) -k + 1}$ in $X$ and use that to show an exponential
growth of the filling of a \text{$(k-1)$-sphere} which lies in a
\text{$k$-flat}, and therefore proving that $\divv_{k-1}(X)$ is
indeed exponential. For more details, see section~$4$
in~\cite{Brady98FIAIFMONC}.

The same idea was taken further by Leuzinger
in~\cite{Leuzinger00CAAFIFSS} to generalize the above result to
any rank~$k$ symmetric space $X$ of nonpositive curvature by
showing the existence of an embedded submanifold $Y \subset X$ of
dimension $n-k+1$ which is quasi-isometric to a Riemannian
manifold with strictly negative sectional curvature, and which
intersects a maximal flat in a geodesic.

Combining Theorem~\ref{T:higherRankSymmetricSpaces} with
Leuzinger's result mentioned above, we obtain the following
corollary which gives a positive answer to
question~{\hyperlink{question.1}{1}}.

\begin{cor}
    The rank of a nonpositively curved symmetric space $X$ can be detected
    using the $\divv_k(X)$ invariants.
\end{cor}

Brady and Farb in~\cite{Brady98FIAIFMONC} suspected that $\divv_1$ has exponential growth for the symmetric space
$SL_n(\mathbb{R})/SO_n(\mathbb{R})$. That was proved by Leuzinger to be true for $n=3$, see~\cite{Leuzinger00CAAFIFSS}. We show that the same result does not hold anymore for $n > 3$. More generally we prove,

\begin{thm} \label{T:div_1ForSymmetricSpaces}
    If $X$ is a symmetric space of nonpositive curvature and rank
    $k \geq 3$, then $\divv_1(X)$ has a quadratic polynomial growth.
\end{thm}

After studying the case of symmetric spaces, we turn our focus to
the class of Hadamard manifolds of pinched negative curvature. We
show,

\begin{thm} \label{T:negativeCurvature}
    If $X$ is a Hadamard manifold with sectional curvature $-b^2 \leq K \leq -a^2 <
    0$, then $\divv_k(X)$ has a polynomial growth of degree at most $k$ for
    every $k \geq 1$.
\end{thm}

A natural question at this point would be whether the same result
holds if we weaken the assumption on the manifold $X$ to be a
rank~$1$ instead of being negatively curved. We give a negative
answer to that question. We give an example of a nonpositively
curved graph manifold of \text{rank~$1$} where $\divv_1$ is
exponential.

All the previously known results about the filling invariants, see~\cite{Brady98FIAIFMONC, Gersten94DI3MG, Leuzinger00CAAFIFSS}, as well as the new results in this paper suggest a connection between these invariants and the connectedness of the Tits boundary of a Hadamard manifold. Theorem~\ref{T:div_1ForSymmetricSpaces} is one example. The non-simply connected Tits boundary of a symmetric space of rank~$2$ is reflected in the the exponential
growth of $\divv_1$, while the quadratic polynomial growth of $\divv_1$ when the rank is bigger than~$2$ is a direct
consequence of the simple connectedness of the Tits boundary. 

Except in a very few special cases,  very little is known in general about which part of the Tits geometry of a nonpositively curved space is preserved under a quasi isometry. Further investigation of the connection between these invariants and the connectedness of the Tits boundary may shed some light on that question. 

The paper will be organized as follows:
Section~\ref{S:Introduction} is an introduction. In
section~\ref{S:Definitions} we give the precise definition of
$\divv_k$. In section~\ref{S:rankOne} the proof of
Theorem~\ref{T:higherRankSymmetricSpaces} in the simple case of rank~$1$ symmetric
spaces is given where the basic idea is illustrated. In
section~\ref{S:higherRankSymmetricSpaces} we prove
Theorem~\ref{T:higherRankSymmetricSpaces} for higher rank
symmetric spaces. The proof of Theorem~\ref{T:div_1ForSymmetricSpaces} is give
in section~\ref{S:firstFillingInvariant}. In section~\ref{S:negativeCurvature} we prove
Theorem~\ref{T:negativeCurvature}. In
section~\ref{S:graphManifolds} we give the graph manifold example
mentioned above. In appendix~\ref{S:DifferentProof} we give a
different and shorter proof of Leuzinger's Theorem, which was
proved in~\cite{Leuzinger00CAAFIFSS}.

\subsection*{Acknowledgments} The author would like to thank his
advisor Chris Croke for many helpful discussions during the
development of this paper, and also like to thank the University
of Bonn, in particular Werner Ballmann for the opportunity to
visit in the summer of 2003 when his interest began in the filling
invariants.

\section{Definitions and Background} \label{S:Definitions}
Let $X^n$ be an \text{$n$-dimensional} Hadamard manifold, by that
we mean a complete simply connected Riemannian manifold with
nonpositive sectional curvature. By Cartan-Hadamard Theorem $X^n$
is diffeomorphic to $\Real^n$. In fact the $\exp_{x_0}$ map at any
point $x_0 \in X$ is a diffeomorphism. For a standard source on
Hadamard manifolds we refer the reader to~\cite{Ballmann85MONC}.

We denote the ideal boundary of $X$ by $X(\infty)$. For any two different
points $p$, $q \in X$, $\overline{pq}$, $\overrightarrow{pq}$ denote respectively
the geodesic segment connecting $p$ to $q$, and the geodesic ray staring at
$p$ and passing through $q$. By $\overrightarrow{pq}(\infty)$ we denote
the limit point in $X(\infty)$ of the ray $\overrightarrow{pq}$.

Let $\overline{X} = X \cup X(\infty)$, if $x_0 \in X$ and $p$, $q \in
\overline{X} \setminus \{ x_0 \}$ then $\angle_{x_0} (p,q)$ is the
angle between the unique geodesic rays connecting $x_0$ to $p$ and
$q$ respectively. If $p$, $q \in X(\infty)$, then $\angle(p,q) =
\sup_{x_0 \in X} \angle_{x_0}(p,q)$ denotes the Tits angle between
$p$ and $q$. If $\angle_{x_0}(p,q) = \angle(p,q) < \pi$ for some
$x_0 \in X$ then the geodesic rays connecting $x_0$ to $p$ and $q$
respectively bound a flat sector. Conversely if those rays bound a
flat sector then $\angle_{x_0}(p,q) = \angle(p,q)$.

Let $S_r(x_0)$, $B_r(x_0)$ and $B^{\circ}_{r}(x_0)$ denote
respectively the distance sphere, the distance ball and the open
distance ball of radius $r$ and center $x_0$. Let $S^k$ and
$B^{k+1}$ denote respectively the unit sphere and the unit ball in
$\Real^{k+1}$.
 Let $C_r(x_0) = X \setminus B^{\circ}_{r}(x_0)$. Projection
along geodesics of $C_r(x_0)$ onto the sphere $S_r(x_0)$ is a
deformation retract, which decreases distances, since the ball
$B_r(x_0)$ is convex and the manifold is nonpositively curved. Any
continuous map $f \colon S^k \To S_r(x_0)$ admits a continuous
extension, ``filling'', $\hat{f} \colon B^{k+1} \To C_r(x_0)$, and
the extension could be chosen to be Lipschitz if $f$ is Lipschitz.

Lipschitz maps are differential almost everywhere. Let
$\abs{D_x(f)}$ denote the Jacobian of $f$ at $x$. The
\text{$k$-volume} of $f$ and the \text{$(k+1)$-volume} of
$\hat{f}$ are defined as follows,
\begin{gather}
    \vol_k(f) = \int_{S^k} \abs{D_x f}, \\
    \vol_{k+1}(\hat{f}) = \int_{B^{k+1}} \abs{D_x \hat{f}}.
\end{gather}

Let $0 < A$ and $0 < \rho \leq 1$ be given. A Lipschitz map $f
\colon S^k \To S_r(x_0)$ is called \text{$A$-admissible} if
$\vol_k(f) \leq A r^k$ and a Lipschitz filling $\hat{f}$ is called
\text{$\rho$-admissible} if $\hat{f}(B^{k+1}) \subset C_{\rho
r}(x_0)$. Let
\begin{equation}
    \delta_{\rho, \, A}^{k} = \sup_{f} \inf_{\hat{f}}
    \vol_{k+1}(\hat{f}),
\end{equation}
where the supremum is taken over all \text{$A$-admissible} maps
and the infimum is taken over all \text{$\rho$-admissible}
fillings.

\begin{defn}
    The invariant $\divv_k(X)$ is
    the two parameter family of functions
    \begin{equation*}
        \divv_k(X) = \{ \delta_{\rho, \, A}^{k} \mid 0 < \rho \leq
        1 \text{ and } 0 < A \}.
    \end{equation*}
\end{defn}

Fix an integer $k \geq 0$, for any two functions $f$, $g \colon
\Real^{+} \To \Real^{+}$, we write $f \preceq_{k} g$ if there
exist two positive constants $a$, $b$ and a polynomial
$p_{k+1}(x)$ of degree at most $k+1$ with a positive leading
coefficient such that $f(x) \leq a g(bx) + p_{k+1}(x)$. Now we
write $f \sim_{k} g$ if $f \preceq_{k} g$ and $g \preceq_{k} f$.
This defines an equivalence relation.

We say that $\divv_k \preceq \divv'_{k}$ if there exist $0 <
\rho_0$, $\rho'_{0} \leq 1$ and $A_0$, $A'_{0} > 0$ such that for
every $\rho < \rho_0$ and $A > A_0$ there exist $\rho' <
\rho'_{0}$ and $A' > A'_{0}$ such that $\delta_{\rho, \, A}^{k}
\preceq_{k} \delta_{\rho, \, A}^{' k}$. We define $\divv_k \sim
\divv'_{k}$ if $\divv_k \preceq \divv'_k$ and $\divv'_k \preceq
\divv_k$. This is an equivalence relation and under this
identification $\divv_k$ is a quasi-isometry invariant
(see~\cite{Brady98FIAIFMONC} for details).

\begin{rem}
    As a quasi-isometry invariant, a polynomial growth rate of $\divv_k$ is only defined up
    to $k+1$. And the reader should view the $k$ polynomial growth
    rate of $\divv_k$ in Theorem~\ref{T:negativeCurvature}
    accordingly.
\end{rem}

The polynomial bound on $\vol_k(f)$ in the definition above is
essential to prevent the possibility of constructing exponential
\text{$k$-volume} maps requiring exponential volume fillings, and
therefore making the filling invariants always exponential. On the
other hand, allowing the filling to be ``slightly'' inside the
ball, i.e.\@ inside $C_{\rho r}(x_0)$ for some fixed $0 < \rho$ is
needed since continuous quasi-isometries map spheres to distorted
spheres.

\begin{rem} \label{R:volumeOfConing}
    Through the paper we will use the following cone construction.
    Given an \text{$n$-dimensional} Hadamard manifold $X$ and a Lipschitz map $f \colon S^{k} \To
    X$, if we cone $f$ from a point $x_0 \in X$ we obtain a Lipschitz extension $\hat{f} \colon D^{k+1} \To
    X$. Using comparison with Euclidean space, it is clear that
    $\vol_{k+1}(\hat{f}) \leq \vol_{k}(f) \sup_{p \in S^k}
    d(x_0,f(p))$. If $X$ has sectional curvature $K \leq -a^2
    <0$, then there exists a constant $C = C(a,n)$ which does not
    depend on $x_0$ nor $f$ such that $\vol_{k+1}(\hat{f}) \leq C
    \vol_{k}(f)$. If $a =1$ this constant can be chosen to be~$1$.
\end{rem}

\section{Rank One Symmetric Spaces} \label{S:rankOne}
The $\divv_0 (X)$ of a Hadamard manifold $X$ is the same as the
``rate of divergence'' of geodesics. Therefore for rank~1
symmetric spaces it is exponential since they have pinched
negative curvature.

In this section we calculate $\divv_k$ where $k \geq 1$ for any
rank one symmetric space. First we start with a lemma which will
be needed in the proof of Theorem~\ref{T:rankOne}.

\begin{lem} \label{L:deformationLemmaInR^n}
    Given a Lipschitz function $f \colon S^k \To \mathbb{R}^n$, where
    $1 \leq k <n$, there exists a constant $c = c(n)$ such that for any
    ball $B_{r}(u_0)$, we can find a map $g \colon S^k \To
    \mathbb{R}^n$ which satisfies the following conditions:
    \begin{enumerate}[\upshape i.]
        \item $f$ is homotopic to $g$.
        \item $\vol_k (g) \leq c \vol_k (f)$.
        \item $g(S^k) \cap B^{\circ}_r (u_0) = \emptyset$.
        \item $f$ and $g$ agree outside $B_{r}(u_0)$.
        \item The \text{$(k+1)$-volume} of the homotopy is bounded by
        $cr \vol_k(f)$.
    \end{enumerate}
\end{lem}
\begin{proof}
    The proof follows the proof of Theorem~$10.3.3$ in
    \cite{Epstein92WPIG}. The idea is to project the part inside the ball $B_r(u_0)$
    to the sphere $S_r(u_0)$. We will project from a point
    within the ball $B_{r/2}(u_0)$, but since projecting from the
    wrong center might increase the volume by a huge factor, we
    average over all possible projections and prove that the
    average is under control, and therefore we have plenty of
    centers from which we can project and still have the volume of
    the new function under control.

    Let $\omega_{n-1}$ be the volume of the unit \text{$(n-1)$-sphere} in
    $\mathbb{R}^n$. Let $|D_x f |$ represent the Jacobian of $f$ at
    the point $x$. Let $\pi_{u}$ be the projection map from the
    point $u \in B_{r/2} (u_0)$ to the sphere $S_r(u_0)$.
    \begin{equation} \label{E:projectionEstimate}
        \vol_k (\pi_u \circ f) \leq \int_{f^{-1}( B_r (u_0))} \dfrac{|D_x f | \, (2r)^k}{\norm{f(x) -
        u}^{k} \cos \theta} \, dx + \vol_{k} (f),
    \end{equation}
    where $\theta$ is the angle between the ray connecting $u$ to
    $f(x)$ and the outward normal vector to the sphere at the point of
    intersection of the ray and the sphere $S_r(u_0)$.
    If $u \in B_{r/2}(u_0)$ then it is easy to see that $\theta \leq
    \pi/6$ and therefore $1/\cos \theta \leq 2/\sqrt{3} \leq 2$.

    By integrating~\eqref{E:projectionEstimate} over the ball
    $B_{r/2}(u_0)$ we get,
    \begin{equation}
    \begin{split}
            \int_{B_{r/2}(u_0)} \vol_k(\pi_u \circ f) \, du \leq& \int_{B_{r/2}(u_0)}
                \int_{f^{-1}(B_r(u_0))} \dfrac{|D_x f | \, (2r)^k}{\| f(x) -u
                \|^k \cos \theta}\, dx \, du\\
            &\qquad + \vol(B_{r/2}(u_0)) \, \vol_k(f)\\
            \leq& \int_{f^{-1}(B_r(u_0))} |D_x f| \int_{B_{r/2}(u_0)}
                \dfrac{2 \, (2r)^k}{\norm{f(x) - u}^k} \, du \, dx\\
            &\qquad{} + \vol(B_{r/2}(u_0)) \, \vol_k(f)\\
            \leq&\ 2^{k+1} r^k \int_{f^{-1}(B_r(u_0))} |D_x f | \int_{B_{3r/2}(f(x))}
                \dfrac{1}{\norm{u}^k} \, du \, dx\\
            &\qquad {}+ \vol(B_{r/2}(u_0)) \, \vol_k(f)\\
            =&\ 2^{k+1} r^k \int_{f^{-1}(B_r(u_0))} |D_x f |
            \int_{S^{n-1}} \int_{0}^{3r/2} \dfrac{1}{t^k} t^{n-1} \, dt
            \, d\mu \, dx\\
            &\qquad {}+ \vol(B_{r/2}(u_0)) \, \vol_k(f)\\
            =&\ \dfrac{2^{k+1} 3^{n-k} r^n \omega_{n-1}}{2^{n-k} (n-k)}  \int_{f^{-1}(B_r(u_0))} |D_x f| \,
            dx\\
            &\qquad{} + \vol(B_{r/2}(u_0)) \, \vol_k(f)\\
            \leq& \left[ \dfrac{2^{2k+1} 3^{n-k} n}{n-k} + 1
            \right] \vol(B_{r/2}(u_0)) \vol_k(f)\\
            \leq& \left[ 2^{2n-1} 3^{n-1} n +1 \right] \vol(B_{r/2}(u_0))
            \vol_k(f).
    \end{split}
    \end{equation}
    The \text{$(k+1)$-volume} of the homotopy is bounded by $2r [ 2^{2n-1} 3^{n-1} n
    +1] \vol_k(f)$. We take $c = 2 [ 2^{2n-1} 3^{n-1} n +1]$, which now satisfies all the
    requirements of the lemma.
\end{proof}

Using the standard notation, we denote by $\sn_k$ respectively
$\cs_k$ the solution to the differential equation $x''(t) + k x(t)
= 0$ with the initial conditions $x(0) = 0$ and $x'(0) =1$
respectively $x(0)=1$ and $x'(0)=0$. We also set $\ct_k =
\cs_k/\sn_k$.

Now we use Lemma~\ref{L:deformationLemmaInR^n} to prove a similar
result for any Riemannian manifold with bounded geometry which we
will call the ``Deformation Lemma''.

\hyperdef{Lemma}{1}{}
\begin{lem}[Deformation Lemma] \label{L:deformationLemma}
    Given a Riemannian manifold M, with sectional curvature $L \leq K \leq
    H$ and injectivity radius $\inj(M) \geq \epsilon$, then there
    exist constants $\delta = \delta(L,K,\epsilon) > 0$ and
    $c=c(n,L,K,\epsilon)$, which
    do not depend on $M$ such that for any Lipschitz map $f \colon S^k \To
    M$ and any
    ball $B_{r}(p) \subset M$ with radius $r \leq \delta$, we can find a map $g \colon S^k \To
    M$ which satisfies the following conditions:
    \begin{enumerate}[\upshape i.]
        \item $f$ is homotopic to $g$.
        \item $\vol_k (g) \leq c \, \vol_k (f)$.
        \item $g(S^k) \cap B^{\circ}_r (p) = \emptyset$.
        \item $f$ and $g$ agree outside $B_r(p)$.
        \item The \text{$(k+1)$-volume} of the homotopy is bounded by $c \,
        \vol_k(f)$.
    \end{enumerate}
\end{lem}
\begin{proof}
    Take $\delta = \min \{ \, \epsilon/2 , \pi/2\sqrt{H} \, \}$. By comparison with the spaces of
    constant curvature $L$ and $H$, we have the following bounds
    within $B_{\delta}(p)$.
    \begin{align}
        \norm{D \exp_{p}^{-1}} &\leq A = \max \{ 1 ,
        \delta/\sn_H(\delta) \}, \label{E:controlOnInverseExpMap} \\
        \norm{D \exp_p} &\leq
        B = \max \{ 1 , \sn_L(\delta)/\delta \}. \label{E:controlOnExpMap}
    \end{align}
    See Theorem~$2.3$ and
    Corollary~$2.4$ in chapter~$6$ of \cite{Petersen98RG} for details but the reader should be aware that the estimate
    for $\norm{D \exp_{p}^{-1}}$ is stated incorrectly there.
    Note that $A$ and $B$ only depend on $L$, $K$ and $\epsilon$ but not
    on $p$ nor $M$.

    Now we use $\exp_p^{-1}$ to lift the map $f$ locally, i.e.\@ in $B_{\delta}(p)$, to $T_p M$
    then use Lemma~\ref{L:deformationLemmaInR^n} to deform this
    lifted map to a map that lands
    outside the ball $B_{\delta}(0) \subset T_{p}M$ and then project back using $\exp_p$ to the
    manifold $M$. The \text{$k$-volume} of the new map and the \text{$(k+1)$-volume}
    of the homotopy will be controlled because of the bounds on  $\norm{D \exp_p}$ and $\norm{D \exp_p^{-1}}$ given
    by equations~\eqref{E:controlOnInverseExpMap} and~\eqref{E:controlOnExpMap} and the
    estimates given in Lemma~\ref{L:deformationLemmaInR^n}. This finishes
    the proof of the lemma.
\end{proof}

\begin{thm} \label{T:rankOne}
    If $X$ is a rank~$1$ symmetric space of noncompact type,
    then $\divv_k(X)$ has a polynomial growth of degree at most $k$ for
    every $k \geq 1$.
\end{thm}

\begin{proof}
    Let $x_0$ be any point in $X$, $S_r(x_0)$ the distance sphere of radius $r$ centered at $x_0$,
    and $d_{S_r(x_0)}$ the Riemannian distance function on the sphere.

    We will show that there exists $0 < \rho < 1$ such that there is a
    filling of any admissible $k$-sphere on $S_r(x_0)$ outside $B^{\circ}_{\rho r}
    (x_0)$ which grows polynomially of degree at most $k$ as $r \to
    \infty$.

    Let $\pi_{r}^{R} \colon S_R(x_0) \To S_r(x_0)$ be the radial
    projection from $S_R(x_0)$ to $S_r(x_0)$, and let
    $\lambda_{r}^{R} \colon S_{r} (x_0) \To S_{R}(x_0)$ be its inverse.
    Assume that the metric on
    $X$ is normalized such that the sectional curvature is bounded
    between $-4$ and $-1$. By comparison with the spaces of constant
    curvatures $-1$ and $-4$, it is easy to
    see that the map $\pi_{r}^{R}$ decreases distance by at least
    a factor of $\sinh R /\sinh r$, and $\lambda_{r}^{R}$
    increases distance at most by a factor of
    $\sinh 2R /\sinh 2r$.

    Fix $A > 0$ and let $f \colon S^k \To S_r(x_0)$ be a Lipschitz map
    such that $\vol_{k}(f) \leq A r^k$. Then $\vol_{k}(\pi_{\rho
    r}^{r} \circ f) \leq A r^k \left( \sinh \rho r /\sinh r
    \right)^k$.

    Since horospheres are Lie groups with left invariant metrics,
    the curvature is bounded above and below and there is a lower
    bound on the injectivity radius. This puts a uniform bound on
    the curvature as well as the injectivity radius of distance
    spheres $S_r(x_0)$ as long as $r$ is big enough. Because of
    that and by using the {\hyperlink{Lemma.1}{Deformation Lemma}},
    we could deform $\pi_{\rho r}^{r} \circ f$ on $S_{\rho r}(x_0)$ to a new function $g$
    which misses a ball, on the sphere $S_{\rho r}(x_0)$, of radius $\delta$
    and has \text{$k$-volume} $ \leq c \vol_{k} (\pi_{\rho r}^{r} \circ f)$
    where $\delta$ and $c$ are the constants given by the
    Deformation Lemma. And such that the homotopy between $\pi_{\rho r}^{r} \circ
    f$ and $g$ has \text{$(k+1)$-volume} $\leq c \vol_{k}(\pi_{\rho r}^{r} \circ
    f)$.

    Let $p \in S_{\rho r} (x_0)$ be the center of this ball and let $q$ be
    its antipodal point. We will cone $g$ from $q$, inside the sphere, and
    then project the
    cone from $x_0$ back to
    the sphere $S_{\rho r}(x_0)$. Because the curvature is less
    than $-1$ the $(k+1)$-volume of the cone is smaller than
    $\vol_{k}(g) \leq c \vol_{k}(\pi_{\rho r}^{r} \circ f)$ (see
    Remark~\ref{R:volumeOfConing}).

    To estimate the volume of the projection of the cone, we need to find
    a lower bound on the distance from $x_0$ to
    the image of the cone. Let $x \in S_{\rho r} (x_0)$ such
    that $d_{S_{\rho r} (x_0)}(x,p) = \delta$. Let $x'$ be the intersection point of
    $S_{2 \rho r} (q)$ and the ray $\overrightarrow{qx}$. For large values of $r$, $x$ and $x'$ are close, and
    therefore we have $d_{S_{2 \rho r} (q)} (p,x') \geq
    \delta/2$. By comparison to the space of constant curvature $-4$ we see
    that $\theta = \angle_{q} (p,x) \geq {\delta}/ (2 \sinh 4 \rho r)$. Now
    by comparison with Euclidian space we have $d(x_0, \overrightarrow{qx}) \geq
    \rho r \sin \theta \geq \rho r \sin(\delta /(2 \sinh 4 \rho
    r))$. So the cone from $q$ misses a ball of radius $\rho r
    \sin (\delta /(2\sinh 4 \rho
    r))$ around $x_0$. Projecting the cone which lies inside $S_{\rho r}(x_0)$ from $x_0$
    to $S_{\rho r}(x_0)$ gives us a
    filling for $g$ of volume $\leq c A r^k \left(
    \sinh \rho r / \sinh r \right)^k \left( \frac{\sinh 2 \rho r}{ \sinh (2
    \rho r \sin (\delta /2 \sinh 4 \rho r))} \right)^{k+1}$.

    Notice that $\sinh t$ behaves like $e^t/2 $ as $t \to \infty $ and
    like $t$ as $t \to 0$, while $\sin t$ behaves like $t$ as $t \to
    0$. Because of that the above estimate of the filling of $g$
    behaves like $c A r^k \left( e^{\rho r}/e^{r} \right)^{k} \left(
    \frac{e^{2 \rho r}/2}{2 \delta \rho r / e^{4 \rho r}}
    \right)^{k+1} =
    \frac{c A}{(4 \delta \rho)^{k+1} r} \, e^{((7k+6)\rho - k)r}$. Clearly we can
    choose $\rho$ small enough to make $(7k+6)\rho - k$ negative.
    Therefore there exists $r_{0} > 0$ such that for all $r \geq
    r_{0}$ the \text{$(k+1)$-volume} of the filling of $g$ will be
    smaller than~$1$. Now the filling of the original map $f$ consists
    of three parts.
    \begin{enumerate}[\upshape i.]
        \item The radial projection of $f$ to the sphere $S_{\rho r}(x_{0})$,
        which has \text{$(k+1)$-volume} $\leq Ar^k$.
        \item The homotopy used to deform the map $\pi_{\rho r}^{r} \circ f$
        to the new map $g$, which has \text{$(k+1)$-volume} $\leq c \vol_{k}(\pi_{\rho r}^{r} \circ f)
        \leq c \vol_{k}(f) \leq cAr^k$.
        \item The filling of $g$, which has \text{$(k+1)$-volume}~$ < 1$.
    \end{enumerate}

    This filling of $f$ is Lipschitz and lies outside $B_{\rho r}^{\circ}(x_0)$, and is therefore
    \text{$\rho$-admissible}, and it has \text{$(k+1)$-volume}
    smaller than $A(c+1)r^k + 1$. This finishes the proof of the theorem.
\end{proof}

\section{Higher Rank Symmetric Spaces} \label{S:higherRankSymmetricSpaces}
In this section we give the proof for
Theorem~\ref{T:higherRankSymmetricSpaces}. The idea is similar to
the one used in the proof of Theorem~\ref{T:rankOne}, namely
project the Lipschitz map you wish to fill, to a smaller sphere to
make the \text{$k$-volume} of the projection small and then fill
the projection on the smaller sphere.

For the rank one symmetric space case, the projection decreases
distances exponentially in all direction, which is no longer true
for the higher rank case. Nevertheless, we still have $n-l$
directions in which the projection, to smaller spheres, decreases
distance by an exponential factor. Hence, the projection will
decrease the volume of a Lipschitz map $f$ by an exponential
factor as long as the dimension of the domain of $f$ is bigger
than $l-1$, to include at least one of the exponentially
decreasing directions. See Lemma~\ref{L:higherRankProjection} for
details.

We start by recalling some basic facts about symmetric spaces. For
a general reference of symmetric spaces of nonpositive curvature
and for the proofs of the facts used in this section see
chapter~$2$ in~\cite{Eberlein96GONCM}.

Let $X = G/K$ be a symmetric space of nonpositive curvature of
dimension $n$ and rank $l$. Fix a point $x_0 \in X$ and let $\gl =
\kl + \pl$ be the corresponding Cartan decomposition, where $\kl$
is the Lie algebra of the isotropy group $K$ at $x_0$. Fix a
maximal abelian subspace $\al \subset \pl$ and let $F$ be the flat
determined by $\al$. We identify $\pl$ with $T_{x_0} X$ in the
usual way. For each regular vector $v \in \al$, let $R_v \colon
v^{\bot} \To v^{\bot}$ denote the curvature tensor where $R_v(w) =
R(w,v)v$. Let $\lambda_1(v), \dots ,\lambda_{n-1}(v)$ be the
eigenvalues of $R_v$ with corresponding eigenvectors $E_1, \dots ,
E_{n-1}$ such that $E_1, \dots ,E_{l-1}$ are tangent to the
$l$-flat $F$. Therefore $\lambda_1(v) = \dots = \lambda_{l-1}(v) =
0$ while $\lambda_{j}(v) < 0$ for $j \geq l$. Since for every
$v_1$, $v_2 \in \al$, $R_{v_1}$ and $R_{v_2}$ commute and
therefore can be simultaneously diagonalized, $E_{l}, \dots ,
E_{n-1} \in (T_{x_0} F)^{\bot}$ can be chosen not to depend on the
choice of $v$ and the $\lambda_j(v)$ are continuous functions in
$v$.

Let $E_j(t)$ be the parallel translates along the geodesic
$\gamma_v (t)$ with $\gamma_v (0) = x_0$ and $\dot{\gamma_v} (0) =
v$. Then $J_j(t) = f_j(t) E_j(t)$ determine an orthogonal basis
for the space of Jacobi fields along $\gamma_v(t)$ which vanishes
at $x_0$ and are orthogonal to $\gamma_v(t)$, where $f_j(t) = t$
if $\lambda_j(v) = 0$ (i.e.\@ $1 \leq j \leq l-1$) and $f_j(t) =
1/ \sqrt{-\lambda_j (v)} \sinh \sqrt{-\lambda_j (v)} t$ if
$\lambda_j(v) < 0$ (i.e.\@ $l \leq j \leq n-1$).

Let $\mathcal{C}$ be the collection of all Weyl chambers of the
first type (i.e.\@ Weyl chambers defined on $S_{x_0} X$), where
$S_{x_0} X \subset T_{x_0}X$ is the sphere of unit vectors in the
tangent space $T_{x_0} X$. We will often identify $S_{x_0} X$ with
$S_r(x_0)$ in an obvious way. Let $\mathscr{C}(v)$ be the Weyl
chamber containing $v$ for any regular vector $v$. Let
$v_{\mathscr{C}}$ be the algebraic centroid of the Weyl chamber
$\mathscr{C}$. Let $\mathcal{B}$ be the set of the algebraic
centroids of all Weyl chambers. By the identification of $S_{x_0}
X$ and $S_r(x_0)$, we will think of $\mathcal{B}$ as a set of
points on $S_r(x_0)$, a set of vectors of length $r$ in $T_{x_0}
X$ and a set of points in $X(\infty)$. The meaning will be clear
from the context. Recall that the subgroup of isometries $K$
fixing $x_0$ acts transitively on $\mathcal{B}$.

Assume $v \in \mathcal{B}$ belongs to the \text{$l$-flat} $F$. Let
$\mathscr{C}(v)$ be the unique Weyl chamber containing $v$. Since
$v$ is a regular vector then $\lambda_l(v), \dots ,
\lambda_{n-1}(v)$ are negative. By the continuity of $\lambda_j$,
there exist $\epsilon > 0$ and $\delta > 0$ such that for every $w
\in N_{\epsilon}(v)$, $\lambda_j(w) \leq - \delta^2$ for every $j
\geq l$, where $N_{\epsilon}(v) = \{ \, w \in \mathscr{C}(v) \mid
\angle_{x_0} (v,w) = \angle(v,w) \leq \epsilon \, \}$. We choose
$\epsilon$ small enough such that $N_{\epsilon}(v)$ is contained
in the interior of $\mathscr{C}(v)$ and away from the boundary of
$\mathscr{C}(v)$. By the transitive action of $K$ on $\mathcal{B}$
we have such neighborhood around each element of $\mathcal{B}$.
Let $\mathcal{N}_{\epsilon}$ be the union of these neighborhoods.
Recall that any two Weyl chambers at infinity are contained in the
boundary at infinity of a flat in $X$. Note that if $v \in
\mathcal{B}$ and $v'$ is the centroid of an opposite Weyl chamber
to $\mathscr{C}(v)$ then $\angle (v, v') = \pi$. And therefore,
for any $v \in \mathcal{B}$ and any $w \notin
\mathcal{N}_{\epsilon}$, we have $\angle (v, w) < \pi - \epsilon$.

From the above discussion we have immediately the following lemma.

\begin{lem} \label{L:higherRankProjection}
    Let $0 < \rho \leq 1$, $k \geq l$ and $\pi_{\rho r}^{r} \colon S_{r}(x_0) \To S_{\rho
    r}(x_0)$ be the projection map. If $f \colon S^k \To S_{r}(x_0)$
    is a Lipschitz map then $\vol_k(\pi_{\rho r}^{r} \circ f \cap \,
    \mathcal{N}_{\epsilon}) \leq \vol_k(f \cap \,
    \mathcal{N}_{\epsilon}) \, (\sinh \rho \delta r / \sinh \delta
    r) \leq \vol_k (f) \, (\sinh \rho \delta r / \sinh \delta r)$.
\end{lem}

Before starting the proof of
Theorem~\ref{T:higherRankSymmetricSpaces} we need the following
lemma.

\begin{lem} \label{L:distanceAngleEstimate}
    Let $X$ be a Hadamard manifold, and let $x_0 \in X$. For every
    $\epsilon > 0$ there exists $\eta > 0$ such that for any two
    points $p_1$, $p_2 \in S_r(x_0)$ with $\angle (\overrightarrow{x_0 p_1}
    (\infty), \overrightarrow{x_0 p_2} (\infty)) \leq \pi - \epsilon$ we have $d(x_0, \overline{p_1 p_2})
    \geq \eta r$.
\end{lem}
\begin{proof}
    Let $\gamma_i = \overrightarrow{x_0 p_i}(t)$ and $z_i = \gamma_i
    (\infty)$ for $i=1$,~$2$. Let $\alpha_i (t) = \angle_{\gamma_i(t)} (x_0,
    \gamma_j(t))$ for $j \not = i$.

    Note that $\pi - (\alpha_1(t) + \alpha_2(t))$ is an increasing
    function of $t$ which converges to $\angle (z_1,z_2)$ as $t \to
    \infty$. Therefore $\epsilon \leq \alpha_1(r) + \alpha_2
    (r)$. Without any loss of generality we may assume that
    $\epsilon/2 \leq \alpha_1(r)$. Using the first law of
    cosines we have $r^2 + d(p_1,p_2)^2 - 2r d(p_1,p_2) \cos (\epsilon/2) \leq
    r^2$. Therefore $d(p_1,p_2) \leq 2 r \cos (\epsilon/2)$. Let
    $m$ be the closest point on $\overline{p_1 p_2}$ to $x_0$, and
    without loss of generality assume that $d(m,p_1) \leq
    d(m,p_2)$. By taking $\eta = 1 - \cos (\epsilon/2)$, we have
    \begin{equation}
        \begin{split}
            d(x_0, \overline{p_1 p_2}) ={}& d(x_0, m)\\
            \geq{}& d(x_0,p_1) - d(p_1,m)\\
            \geq{}& r - r \cos (\epsilon/2)\\
            ={}& \eta r,
        \end{split}
    \end{equation}
    which finishes the proof of the lemma.
\end{proof}

\begin{proof}[Proof of Theorem~\ref{T:higherRankSymmetricSpaces}]
    Let $- \lambda^2$ be a lower bound on the sectional curvature of $X$.
    Fix $A > 0$ and let $f \colon S^k \To S_r(x_0)$ be a Lipschitz map
    such that $\vol_{k}(f) \leq A r^k$. By
    Lemma~\ref{L:higherRankProjection}, $\vol_k(\pi_{\rho r}^{r} \circ f \cap
    \, \mathcal{N}_{\epsilon}) \leq \vol_k (f) \, (\sinh \rho \delta r / \sinh \delta
    r)$.

    Fix $q \in S_{\rho r}(x_0)$ to be the algebraic centroid of a Weyl chamber, and
    let $p \in S_{\rho r}(x_0)$ be its antipodal point. Since curvature and injectivity
    radii of large spheres are again controlled, using the {\hyperlink{Lemma.1}{Deformation Lemma}} as in
    the proof of Theorem~\ref{T:rankOne}, we can assume that $\pi_{\rho r}^{r} \circ
    f$ misses a ball $B_{\mu}(p)$ on $S_{\rho r}(x_0)$ and $\mu$
    does not depends on $r$ as long as $r$ is large enough.

    Now cone the image of $\pi_{\rho r}^{r} \circ f$ from $q$ inside the sphere
    $S_{\rho r}(x_0)$. Let $C_f$ be the image of the cone.
    Let $C_1$ be the part
    of the cone coming from $\pi_{\rho r}^{r} \circ f \cap \,
    \mathcal{N}_{\epsilon}$ and $C_2 = C_f \setminus C_1$. By
    Lemma~\ref{L:higherRankProjection}, $\vol_{k+1}(C_1) \leq 2 \rho Ar^{k+1}
    (\sinh \rho \delta r / \sinh \delta r) \leq 2 Ar^{k+1}
    (\sinh \rho \delta r / \sinh \delta r) $ and $\vol_{k+1} (C_2) \leq
    2\rho Ar^{k+1} \leq 2 Ar^{k+1}$. Taking $\eta$ to be the constant given by
    Lemma~\ref{L:distanceAngleEstimate}, we see that
    $C_2 \cap B_{\eta \rho r}^{\circ}(x_0) = \emptyset$.

    By comparison with the space of constant curvature
    $-\lambda^2$ and arguing as in the proof of Theorem~\ref{T:rankOne},
    it is easy to see that $C_1$ misses a ball of
    radius $\rho r \sin (\mu /(2 \sinh 2 \lambda \rho
    r))$ around $x_0$. Now we project from $x_0$ the part of $C_1$ in $B_{\eta \rho r}(x_0)$
    to the sphere $S_{\eta \rho r}(x_0)$ to obtain a filling $g$ of $\pi_{\rho r}^{r} \circ f$
    lying outside $B_{\eta \rho r}^{\circ}(x_0)$. Now
    \begin{equation}
        \begin{split}
            \vol_{k+1}(g) \leq&\ \vol_{k+1} (C_2) + \vol_{k+1}(C_1) \left( \dfrac{\sinh \lambda \eta \rho
            r}{\sinh (\lambda \rho r \sin( \mu /(2 \sinh 2 \lambda \rho
            r)))} \right)^{k+1}\\
            \leq&\ 2Ar^{k+1} \left[ 1 +
            \left( \dfrac{\sinh \rho \delta r}{\sinh \delta r} \right)
            \, \left( \dfrac{\sinh \lambda \eta \rho
            r}{\sinh (\lambda \rho r \sin(\mu /(2 \sinh 2 \lambda \rho
            r)))} \right)^{k+1} \right].
        \end{split}
    \end{equation}

    This estimate behaves like
    $2Ar^{k+1} [1 + e^{(\rho( \delta + \lambda (k+1)(2+ \eta)) - \delta )
    r}/(2\lambda \rho \mu)^{k+1} r^{k+1}]$ as $r \to \infty$. Clearly we can choose
    $\rho$ small enough to make the exponent $\rho( \delta + \lambda (k+1)(2+ \eta)) -
    \delta$ negative and therefore the \text{$(k+1)$-volume} of the \text{$\eta \rho$-admissible}
    filling of $\pi_{\rho r}^{r} \circ
    f$ will be bounded by $4Ar^{k+1}$ and the \text{(k+1)-volume} of the \text{$\eta \rho$-admissible}
    filling of $f$ is bounded by $5Ar^{k+1}$. This finishes the
    proof of the theorem where $\eta \rho$ is taken in place of $\rho$.
\end{proof}

\section{First filling Invariant for Symmetric Spaces}
\label{S:firstFillingInvariant}
In this section we prove Theorem~\ref{T:div_1ForSymmetricSpaces}. The proof consists
of two steps.

The first step is to deform the closed curve we wish to
fill to a new curve which is continuous when viewed as a curve at
infinity with respect to the Tits metric. The length of the new
curve and the area of the deformation will have to be under
control. The proof is valid for any symmetric space of higher
rank. We will establish this in
Lemma~\ref{L:DeformationInSymmetricSpaces}.

In the second step we will use the assumption that the rank is
bigger than or equal $3$, to fill the new curve with a disk with
controlled area. This step will be done in
Proposition~\ref{P:isopermetricInequaltyForBuildings}.

We denote by $\lambda_{r}^{\infty} \colon S_{r}(x_0)
\longrightarrow (X(\infty), Td)$ the map obtained by sending a
point $x \in S_r(x_0)$ to the point $\overrightarrow{x_0 x}
(\infty) \in X(\infty)$. Notice that this map is almost never
continuous.

\begin{lem} \label{L:DeformationInSymmetricSpaces}
    Let $X$ be a symmetric space of noncompact type and rank $k \geq
    2$. There exist two constants $c > 0$ and $0 < \rho \leq 1$
    which depend only on $X$ such that, for every
    Lipschitz curve $f \colon S^1 \longrightarrow S_{r}(x_0)$, there
    exists a new curve $g \colon S^1 \longrightarrow S_{r}(x_0)$
    which satisfies the following conditions:
        \begin{enumerate}[\upshape i.]
        \item $\lambda_{r}^{\infty} \circ g$ is continuous with
              respect to the Tits metric.
        \item There exists a homotopy between $f$ and $g$ which lies outside $B_{\rho
        r}^{\circ}(x_0)$.
        \item The area of the homotopy is bounded above by $cr \vol_1(f)$.
        \item $\vol_1 (g) \leq c \vol_1 (f)$.
    \end{enumerate}
\end{lem}
\begin{proof}
    Let $\delta = \min (1, \eta/2)$, where $\eta$ is the constant given by
    Lemma~\ref{L:distanceAngleEstimate} for $\epsilon = \pi/2$.
    Divide $f$ into pieces each of
    length $\delta r$, the last piece possibly could be shorter than $\delta r$. If that is
    the case we will call it ``short'' and all the other pieces ``long''.
    
    Let $c_j \colon [0,1] \longrightarrow
    S_r(x_0)$ be one of these pieces. Let $\mathscr{C}_i$ be a
    Weyl chamber containing $c_j(i)$ for $i=0, 1$. Let
    $\mathscr{A}$ be an apartment containing $\mathscr{C}_0$ and
    $\mathscr{C}_1$. Let $\gamma_j \colon [0,1] \longrightarrow S_r (x_0)$
    be a geodesic (with respect to the Tits metric)
    in the apartment $\mathscr{A}$ connecting $c_j(0)$ to $c_j(1)$.
    The goal is to replace the piece $c_j$ with the geodesic
    $\gamma_j$. The homotopy between them which leaves the end
    points $c_j(0)$ and $c_j(1)$ fixed will be through geodesics in
    $X$ connecting $c_j(t)$ to $\gamma_j(t)$ for every $0 \leq t \leq
    1$. The length of these geodesics is no longer than $2r$, since they lie 
    inside the ball $B_r(x_0)$.

    The new curve $g$ will be formed by replacing each piece $c_j$
    with the curve $\gamma_j$. Notice that the length of each $\gamma_j$
    is no bigger than $\pi r$. If $\vol_1(f) < \delta r$, i.e.\@ we have no ``long''
    pieces then $g$ is just a 
    point and the statement trivially follows. If $\vol_1(f) \geq \delta r$, i.e.\@ we have at 
    least one long piece then
    it is not hard to see that $\vol_1(g) \leq 2 \pi \vol_1(f)/\delta$.

    Since $X$ is nonpositively curved, the
    area of the homotopy between $c_j$ and $\gamma_j$ is no bigger than $2 \pi r^2$. And
    the area of the homotopy between $f$ and $g$ is no bigger than $cr \vol_1(f)$, where
    $c=4\pi/\delta$.

    To finish the proof we need to show that the homotopy lies
    outside $B_{\rho r}^{\circ}(x_0)$. So we need to show that the
    distance between $x_0$ and the
    geodesic, in $X$, connecting $c_j(t)$ and $\gamma_j(t)$ is at
    least $\rho r$ for every $0 \leq t \leq 1$. Assume that $t \leq
    1/2$, the other case is similar. Notice that $\angle(c_j(0), \gamma_j(t)) \leq
    \pi/2$. Applying Lemma~\ref{L:distanceAngleEstimate}, we get
    that $d_X(x_0, \overline{c_j(0) \gamma_j(t)}) \geq \eta r$,
    where $\eta$ is the constant in Lemma~\ref{L:distanceAngleEstimate}.

    Recall that $d_X(c_j(0), c_j(t)) \leq \delta  r$. For every point $s$ on the
    geodesic $\overline{c_j(t) \gamma_j(t)}$ there exists a point $s'$ on the geodesic
    $\overline{\gamma_j(t) c_j(0)}$ such that $d_X(s,s') \leq \delta
    r$. Therefore $d_X(x_0,s) \geq (\eta - \delta)r \geq \eta r/2$. By taking
    $\rho = \eta/2$, the image of the homotopy lies
    outside $B_{\rho r}^{\circ}(x_0)$. This finishes the proof of
    the lemma.
\end{proof}

Now we proceed to the second step of the proof. We prove a more
general result.

\begin{prop} \label{P:isopermetricInequaltyForBuildings}
    Let $\Delta$ be a spherical building with a re-scaled metric such that
    each apartment is isometric to $S^{n-1}(r)$, the round sphere of radius $r$.
    There exists a constant $c > 0$
    which depends only on $\Delta$ but not $r$ such that for any Lipschitz 
    function $g \colon S^{k} \longrightarrow \Delta$, where $k < n-1$,
    we can extend $g$ to a new function $\hat{g} \colon B^{k+1} \longrightarrow
    \Delta$ such that $\vol_{k+1}(\hat{g}) \leq cr \vol_{k}(g)$.
\end{prop}
\begin{proof}
    Fix a point $p$ to be the center of a Weyl chamber
    $\mathscr{C}$. Let $\mathcal{A}$ be the collection of all opposite Weyl
    chambers to $\mathscr{C}$, and let $\mathcal{B}$ be the collection of
    all antipodal points to $p$. Let $\epsilon > 0$ be the largest positive number
    such for any $q \in \mathcal{B}$ the ball $B_{\epsilon r}(q)$ is contained
    in the interior of the Weyl chamber containing $q$. Notice that $\epsilon$
    only depends on $\Delta$.

    The proof of Lemma~\ref{L:deformationLemmaInR^n} can be easily modified
    to deform the function $g$ to miss the ball $B_{\epsilon r}(q)$ for every $q \in \mathcal{B}$.
    Now we cone the new deformed function from the point $p$ to
    obtain the desired extension.
\end{proof}

\begin{proof}[Proof of Theorem~\ref{T:div_1ForSymmetricSpaces}]
The proof is immediate. We deform the function $f$ we wish to fill
to a new function $g$ using
Lemma~\ref{L:DeformationInSymmetricSpaces}. We fill $g$ with a
disc by invoking Proposition~\ref{P:isopermetricInequaltyForBuildings},
for $k = 1$, where we identify $\Delta$ with the Tits building structure on the sphere
$S_r(x_0)$ induced from $(X(\infty), Td)$. The filling of $g$ lies on $S_r(x_0)$ and therefore
outside $B_r^{\circ}(x_0)$.
\end{proof}

\begin{rem}
    We expect a similar result to
    Theorem~\ref{T:div_1ForSymmetricSpaces} to hold for a larger
    class than Symmetric spaces. One candidate is the class of Hadamard
    manifolds whose boundary at infinity is simply connected with
    respect to the Tits metric.
\end{rem}

\section{Riemannian Manifolds of Pinched Negative curvature}
\label{S:negativeCurvature}
In this section we give the proof of
Theorem~\ref{T:negativeCurvature}. All the steps of the proof of
Theorem~\ref{T:rankOne} carry over automatically to our new
setting, except the uniform bounds on the curvature and the
injectivity radius of distance spheres $S_r(x_0)$. By uniform we
mean independent of $r$ for large values of $r$. These bounds will
be established in Lemma~\ref{L:boundsOnCurvature} and
Lemma~\ref{L:boundOnInjectivityRadius} below. The proofs are
straightforward and we include them for completeness. First we
recall proposition~$2.5$ in chapter~IV from
\cite{Ballmann95LOSONC} concerning estimates on Jacobi fields.

\begin{prop}
    Let $\gamma \colon \Real \To M$ be a unit speed geodesic and
    suppose that the sectional curvature of $M$ along $\gamma$ is
    bounded from below by a constant $\lambda$. If $J$ is a Jacobi
    field along $\gamma$ with $J(0) = 0$, $J'(0) \bot
    \dot{\gamma}(0)$ and $\norm{J'(0)} = 1$, then
    \[
        \norm{J(t)} \leq \sn_{\lambda}(t) \text{ and } \norm{J'(t)} \leq
        \ct_{\lambda}(t) \norm{J(t)},
    \]
    if there is no pair of conjugate points along $\gamma |
    [0,t]$.
\end{prop}

    Notice that the estimate on $\norm{J'(t)}$ is still valid without
    requiring that $\norm{J'(0)} = 1$. From the proposition we
    have the following immediate corollary.

\begin{cor} \label{C:boundOnJacobiFields}
    Let $X$ be a Hadamard manifold, with sectional curvature $-b^2 \leq K \leq -a^2 <
    0$, if $\gamma$ is a unit speed geodesic, $J(t)$ is a Jacobi
    field along $\gamma$ with $J(0) = 0$, $J(t) \bot
    \dot{\gamma}(t)$ and $\norm{J(r)} = 1$ then $\norm{J'(r)} \leq \ct_{-b^2}
    (r)$, and therefore the upper bound goes uniformly, independent of $\gamma$,
    to $b$ as $r \to \infty$.
\end{cor}

We use the estimate on the derivative of Jacobi fields to put an
estimate on the second fundamental form of distance spheres, and
hence bounds on the sectional curvature.

\begin{lem} \label{L:boundsOnCurvature}
    Let $X$ be a Hadamard manifold, with sectional curvature $-b^2 \leq K \leq -a^2 <
    0$, and $x_0 \in X$. There
    exists a number $H = H(a,b) > 0$ such that the absolute value of the
    sectional curvature of $S_r(x_0)$ is bounded by $H$ for every
    $r \geq 1$.
\end{lem}
\begin{proof}
    The lemma is immediate using the Gauss formula for the sectional curvature of
    hypersurfaces, namely
    \[
        K(Y,Z) - \overline{K}(Y,Z) = \langle B(Y,Y), B(Z,Z)
        \rangle - \norm{B(Y,Z)}^2.
    \]
    Where $Y$ and $Z$ are orthogonal unit vectors tangent to the hypersurface.

    Notice that for $S_r(x_0)$ we have $\norm{B(Y,Z)} \leq \norm {J'_{Y}(r)}$,
    with $J(0) = 0$ and $J_{Y}(r) = Y$, for every two unit vectors $Y$ and $Z$ tangent to the sphere.
    But $\norm{J'_{Y}(r)}$ is bounded by $b$ for large values of $r$
    by Corollary~\ref{C:boundOnJacobiFields}.
\end{proof}

In the next lemma we establish a lower bound on the injectivity
radius $\inj$ of distance spheres.

\begin{lem} \label{L:boundOnInjectivityRadius}
    Let $X$ be a Hadamard manifold, with sectional curvature $-b^2 \leq K \leq -a^2 <
    0$, and $x_0 \in X$. There exists a number $\delta >
    0$ such that $\inj(S_r(x_0)) \geq \delta$ for every $r \geq 1$.
\end{lem}
\begin{proof}
    Using Lemma~\ref{L:boundsOnCurvature} we
    have a uniform upper bound $H$ on the sectional curvature of
    distance spheres $S_r(x_0)$ for $r \geq 1$ and therefore a  lower bound $\pi/\sqrt{H}$
    on the conjugate radius.

    Our plan is to choose $\delta$ small enough and prove that if $\inj(S_r(x_0)) \geq
    \delta$ then $\inj(S_{r+s}(x_0)) \geq \delta$ for all $0 \leq s \leq
    1$. Using the lower bound $-b^2$ on the sectional curvature of $X$ and by comparison with
    the space of constant curvature $-b^2$, there
    exists a constant $0 < B = B(b)  < 1$ independent of $r$ and $0 \leq s \leq 1$
    such that $\pi_{r}^{r+s} \colon S_{r+s}(x_0) \To S_{r}
    (x_0)$ satisfies the following inequality
    \begin{equation} \label{E:quasiInequalty}
        B \, d(x,y) \leq d(\pi_{r}^{r+s}(x), \pi_{r}^{r+s}(y)) \leq
        d(x,y), \quad \forall \, x, y \in S_{r+s}(x_0).
    \end{equation}
    Assume $\delta < \min\{ \, \inj(S_1(x_0)), B\pi/2\sqrt{H} \, \}$.
    Let us assume that $\inj(S_r(x_0)) \geq \delta $ but $\inj(S_{r+s}(x_0)) <
    \delta$. Since $\delta \leq \pi/2\sqrt{H}$, there exists two
    points $p$, $q \in S_{r+s}(x_0)$ and two minimizing geodesic
    (with respect to the induced Riemannian metric on
    $S_{r+s}(x_0)$) $\gamma_1$, $\gamma_2$ connecting $p$ to $q$.
    Moreover $d(p,q) < \delta$. By Klingenberg's Lemma, see~\cite{doCarmo92RG}, any
    homotopy from $\gamma_1$ to $\gamma_2$ with the end points fixed must contain a curve
    which goes outside $B_{\pi/\sqrt{H}}^{\circ}(p)$. Therefore
    these two curves are not homotopic inside the ball $B_{\delta
    /B}^{\circ} (p)$. Using \eqref{E:quasiInequalty} we have,
    \begin{equation}
        B_{\delta}^{\circ}(\pi_{r}^{r+s}(p)) \subseteq
        \pi_{r}^{r+s}(B_{\delta /B}^{\circ} (p)).
    \end{equation}
    Notice that $\pi_{r}^{r+s} \circ \gamma_i$ are contained
    in $B_{\delta }^{\circ}(\pi_{r}^{r+s}(p))$, but not homotopic to each
    other within that ball, which
    is a topological ball since $\inj(S_r(x_0)) \geq \delta$ and this is a contradiction.
    This finishes the proof of the lemma.
\end{proof}

\section{The Graph Manifold Example} \label{S:graphManifolds}
In this section we give the example mentioned in the introduction
showing that Theorem~\ref{T:negativeCurvature} is false if we
relaxed the condition on the manifold from being negatively curved
to being merely rank~$1$.

Our example will be a graph manifold. Graph manifolds of
nonpositive curvature form an interesting class of
\text{$3$-dimensional} manifolds for various reasons. They are the
easiest nontrivial examples of rank~$1$ manifolds whose
fundamental group is not hyperbolic. They are rank~$1$, yet still
have a lot of \text{$0$-curvature}. In fact every tangent vector
$v \in T_p X$ is contained in a \text{$2$-plane} $\sigma \subset
T_p X$ with curvature $K(\sigma) = 0$.

They were used by Gromov in~\cite{Gromov78MONC} to give examples
of open manifolds with curvature $- a^2 \leq K \leq 0$ and finite
volume which have infinite topological type, contrary to the case
of pinched negative curvature. The compact ones have been
extensively studied by Schroeder in~\cite{Schroeder86RONCG}.

Since we are mainly interested in giving a counter example, and
for the sake of clarity, we will consider the simplest possible
graph manifold. Even though the same idea works for a large class
of graph manifolds.

We start by giving a description of the manifold. Let $W_1$ and
$W_2$ be two tori with one disk removed from each one of them. Let
$ B_i = W_i \times S^1$. Each $B_i$ is called a block. The
boundary of each block is diffeomorphic to $S^1 \times S^1$. The
manifold $M$ is obtained by gluing the two blocks $B_1$ and $B_2$
along the boundary after interchanging the \text{$S^1$-factors}.

Since the invariants $\divv_k$ do not depend on the metric, we
choose a metric which is convenient to work with. Take the flat
torus, corresponding to the lattice $\mathbb{Z} \times \mathbb{Z}
\subset \Real^{2}$. Let $\beta_1$ and $\beta_2$ be the unique
closed geodesics of length~$1$. We picture the torus as the unit
square $[-\frac{1}{2}, \frac{1}{2}] \times [-\frac{1}{2},
\frac{1}{2}] \subset \Real^2$ with the boundary identified. Remove
a small disk in the middle and pull the boundary up, such that the
metric is rotationally symmetric around the \text{$z$-axis}. It is
easy to see we can obtain a metric on $W_{i}$ with curvature $
-b^2 \leq K \leq 0$, and make it product near the boundary which
is a closed geodesic. Rescale the metric on $W_i$ to make the
curvature $-1 \leq K \leq 0$. The closed curves $\beta_1$ and
$\beta_2$ are still closed geodesics, and generate the fundamental
group of the punctured torus. Take the metric on each block to be
the product of this metric with a circle of length equal to the
length of the boundary of $W_i$. Gluing the two metrics on the two
blocks together gives a smooth metric on $M$ with curvature $-1
\leq K \leq 0$.

Let $X$ be the universal covering space of $M$, and $\pi \colon X
\To M$ be the covering map. We show that $\divv_1(X)$ is
exponential.

Let $Y$ be a connected component of $\pi^{-1}(B_1)$. Since $B_1$
is a convex subset of $M$, it is easy to see that $Y = Z \times
\Real$ is the universal covering space of $B_1$, where $Z$ is the
universal covering space of $W_1$, and the restriction of $\pi$ to
$Z\times \{0\}$ is the covering map onto $W_1$. We will identify
$Z\times \{0\}$ with $Z$. Clearly $\pi_1(B_1) = \pi_1(W_1) \times
\mathbb{Z}$ and $\pi_{1}(W_1) = \mathbb{Z} \ast \mathbb{Z} $ is a
free group generated by the closed geodesics $\beta_1$ and
$\beta_2$. Moreover the universal covering space $Z$ is a
thickening of the Cayley graph of $\mathbb{Z} \ast \mathbb{Z}$, and it retracts to
it.

Let $w_0$ be the unique intersection point of $\beta_1$ and
$\beta_2$. Take $w_0$ to be the base point of the fundamental
group of $W_1$. Let $a_1$, $a_2 \in \pi_1(W_1, w_0)$ represent the
elements corresponding to $\beta_1$ and $\beta_2$.

Let $\psi \colon Z \To W_1$ be the covering map which is the
restriction to $Z = Z\times\{0\}$ of the covering map $\pi \colon
X \To M$. Fix $p_0 \in Z$ such that $\psi(p_0) = w_0$. We denote
the deck transformation corresponding to any element $s \in
\pi_1(W_1, w_0)$ by $\phi_{s}$.

Let $\gamma_i$ be the lift of $\beta_i$ to a geodesic in $Z$
starting at $p_0$. Notice that $\phi_{a_i}$ is a translation along
the geodesic $\gamma_i$. Take $F = \gamma_1 \times \Real$, which
is a totally geodesic submanifold isometric to $\Real^2$. And let
$x_0 = (p_0, 0) \in Z \times \Real$, where $p_0 = \gamma_1(0) =
\gamma_2(0)$. Let $S_F(r) = S_r(x_0) \cap F$, be the
\text{$1$-sphere} in F with center $x_0$ and radius $r$, $A_Z(r) =
S_r(x_0) \cap Z$ and $B_Z(r) =  B_r(x_0) \cap Z$.

Let $f_r \colon S^1 \To S_F(r)$ be the canonical diffeomorphism
with constant velocity. Notice that $f_r$ is a Lipschitz map and
$\vol_1(f_r) = 2\pi r$ and therefore an admissible map.

Our goal is to show that for every fixed $0 < \rho \leq 1$ the
infimum over all \text{$\rho$-admissible} fillings of $f_r$ grows
exponentially as $r \to \infty$. We show this first for $\rho =1$
and then the general case will follow easily from that.

Notice that since the curvature is nonpositive then any filling of $f_r$ outside $B_r^{\circ}(x_0)$
can be made smaller by radial projection. Therefore
a smallest filling of $f_r$ outside $B_r^{\circ}(x_0)$
would actually lie on the sphere $S_r(x_0)$, so we will only look
at those fillings which lie on the sphere. $S_r(x_0)$ is a
\text{$2$-dimensional} sphere and $f_r$ is a simple closed curve
on it, therefore it divides the sphere into two halves $H_{1}S$
and $H_{2}S$. Any filling of $f_r$ has to cover one of these two
halves. So it is enough to show that $\vol_2(H_{i}S)$ grows
exponentially with $r$ for $i=1$, $2$.

We will estimate the area of each half from below by estimating
the area of the part which lies inside $Y$. The geodesic
$\gamma_1$ divides $Z$ into two halves $H_1$ and $H_2$. We
concentrate on one of them, say $H_1$. Let $b(r) = \vol_2(B_Z(r)
\cap H_1)$, which is an increasing function.

The portion of the half sphere $H_{1}S$ inside $Y$ has area bigger
than the area  of $B_Z(r) \cap H_1$. That is easy to see since
vertical projection of that part will cover $B_Z(r) \cap H_1$, and
the projection from $Y = Z\times \Real$ onto $Z$ decreases
distance since the metric is a product. So it is enough to show
that $b(r) = \vol_2(B_Z(r) \cap H_1)$ grows exponentially.

It is easy to see that $\vol_{2}(B_Z(r))$ grows exponentially with
$r$. One way to see it is to consider the covering map $\psi
\colon Z \To W_1$, and look at all lifts, under deck
transformations, of a small ball around $w_0 \in W_1$. The number
of these disjoint lifted balls which are contained in $B_Z(r)$
increases exponentially because the fundamental group is free and
therefore has exponential growth.

We need to show that the number of these lifted balls in each half
$H_1$ or $H_2$ increases exponentially. We show that for $H_1$.

Without loss of generality we assume that $\phi_{a_2}(p_0) \in
H_1$. We claim that under the deck transformations corresponding
to the subset $S = \{ sa_2 \mid sa_2 \text{ is a reduced word} \}
\subset \pi_{1}(W_1, w_0)$, the image of $p_0$ is in $H_1$. Recall
that $Z$ is a thickening of the the Cayley graph of $\mathbb{Z}\ast \mathbb{Z}$ and
the action of the deck transformation corresponding to the action
of the free group on its Cayley graph. Now the statement follows
since the Cayley graph is a tree. But the number of elements in
$S$ of length less than or equal $m$ increases exponentially as $m
\to \infty$. This finishes the proof of the claim. Therefore
$\vol_2(H_1 S) \geq b(r)$ grows exponentially and that finishes
the proof for $\rho =1$.

Now we turn to the general case where $0 < \rho \leq 1$. We showed
that there exists $0 < \epsilon$ such that for large values of
$r$, $e^{\epsilon r} \leq \inf( \vol_2(\hat{f_r}))$, where the
infimum is taken over all \text{$1$-admissible} fillings, i.e.\@
the fillings which lie in $C_r(x_0) = X \setminus
B_{r}^{\circ}(x_0)$. Fix $0 < \rho \leq 1$ and let $\pi_{\rho
r}^{r} \colon S_r(x_0) \To S_{\rho r}(x_0)$, be the projection
map. Let $g$ be any \text{$\rho$-admissible} filling of $f_r$.
Clearly $\pi_{\rho r}^{r} \circ g$ is a \text{$1$-admissible}
filling of $\pi_{\rho r}^{r} \circ f_r = f_{\rho r}$. Therefore
$\vol_2(g) \geq \vol_2(\pi_{\rho r}^{r} \circ g) \geq e^{\epsilon
\rho r}$, as long as $r$ is big enough. And this finishes the
proof.

\begin{rem}
    In~\cite{Gersten94DI3MG} Gersten studied the growth rate of
    the $\divv_0$ invariant for a
    large class of \text{$3$-manifolds} including graph manifolds.
    Gersten showed in Theorem~5 that a closed Haken
    \text{$3$-manifold} is a graph manifold if and only if the
    $\divv_0$ invariant has a quadratic growth. The author would
    like to thank Bruce Kleiner for bringing this paper to his
    attention.
\end{rem}

\appendix
\section{} \label{S:DifferentProof}

In this appendix we give a different proof of the following
theorem.

\begin{thm}[Leuzinger~\cite{Leuzinger00CAAFIFSS}] \label{T:ExponentialDivForSymmetricSpaces}
    If $X$ is a rank~$k$ symmetric space of nonpositive curvature,
    then $\divv_{k-1}(X)$ has exponential growth.
\end{thm}
\begin{proof}
We use the same notation as in
section~\ref{S:higherRankSymmetricSpaces}. Fix a maximal
\text{$k$-flat} $F$ passing through $x_0$. Fix a Weyl chamber
$\mathscr{C}$ in $F$. Let $v \in \mathscr{C}$ be the algebraic
centroid of $\mathscr{C}$. We identify $X(\infty)$ and $S_r(x_0)$.
Since the set of regular points is an open subset of $S_r(x_0)$
with the cone topology, we can find $0 < \eta$ such that the set S
= $\{ w \in S_r(x_0) \mid \angle_{x_0} (v, w) < \eta \}$ does not
contain any singular point. Moreover by choosing $\eta$ small
enough we may assume that $S \subset \mathcal{N}_{\epsilon}$. To
see this it is enough to show a small neighborhood (with respect
to the cone topology) of $v$ is contained in
$\mathcal{N}_{\epsilon}$. Take any sequence $\{ v_n \}$ which
converges to $v$ in the cone topology. We may assume that $v_n$ is
a regular point for each $n$ since the set of regular points is
open. Any open Weyl chamber is a fundamental domain of the action
of $K$ on the set of regular points in $X(\infty)$, see
proposition~$2.17.24$ in~\cite{Eberlein96GONCM}. Note that the
algebraic centroid of a Weyl chamber is mapped to the algebraic
centroid of another Weyl chamber under the action of $K$.
Therefore for large values of $n$, $v_n$ would be close to the
algebraic centroid of the unique Weyl chamber containing $v_n$,
and therefore contained in $\mathcal{N}_{\epsilon}$.

Let $f_r \colon S^{k-1} \To S_r(x_0) \cap F$ be a diffeomorphism.
We will show that any filling of $f_r$ grows exponentially with $r
\to \infty$. Assume not, then for each $n \in \mathbb{N}$ there is
a filling $\hat{f_n}$ for $f_n$ such that the $\vol_k(\hat{f_n})$
grows sub-exponentially. Let $\pi_{1}^{n} \colon S_{n}(x_0) \To
S_1(x_0)$ be the projection map. By
Lemma~\ref{L:higherRankProjection} $\vol_{k}(\pi_{1}^{n} \circ
\hat{f_n} \cap \mathcal{N}_{\epsilon}) \leq \vol_{k}(\hat{f_n})
(\sinh \delta / \sinh \delta n)$, which goes to zero as $n \to
\infty$. Let $\phi$ be the projection to the flat~$F$. And $g_n$
be the projection of $\pi_{1}^{n} \circ \hat{f_n}$ to $F$. The
image of $g_n$ has to cover the closed unit ball in $F$.

For every $w \notin \mathcal{N}_{\epsilon}$, $\angle_{x_0} (v,w)
\geq \eta$ and therefore $d(\phi(w), v) \geq 2a = 1- \cos \eta$.
Let $A = B_{a}(v) \cap B_1(x_0)$ be the part of the
\text{$k$-ball} in the flat $F$ centered at $v$ with radius $a$
which lies inside the unit closed \text{$k$-ball} in $F$. Now it
is easy to see that the only part of $g_n$ which will lie inside
$A$ would be coming from the portion inside
$\mathcal{N}_{\epsilon}$. The \text{$k$-volume} of that part goes
to zero as $n \to \infty$, nevertheless it has to cover $A$ which
is a contradiction. Therefore the filling of $f_r$ grows
exponentially.
\end{proof}

\bibliographystyle{amsalpha}
\bibliography{bibwith}
\end{document}